\documentclass[a4paper,12pt]{article}
\usepackage{amscd,amsfonts,amssymb,amscd,amsmath,latexsym}
\usepackage[hypertex,backref]{hyperref}
\usepackage[all]{xy}
\begin{titlepage}
\author{Ulrich Bunke}
\title{Orbifold index and equivariant $K$-homology}

\end{titlepage}

\begin{document}

\maketitle

\newcommand{\diag}{{\rm diag}}
\newcommand{\proof}{{\it Proof.$\:\:\:\:$}}
 \newcommand{\dist}{{\rm dist}}
\newcommand{\kaaa}{{\frak k}}
\newcommand{\paaa}{{\frak p}}
\newcommand{\vp}{{\varphi}}
\newcommand{\taaa}{{\frak t}}
\newcommand{\haaa}{{\frak h}}
\newcommand{\R}{{\mathbb R}}
\newcommand{\Hh}{{\bf H}}
\newcommand{\Q}{{\mathbb{Q}}}
\newcommand{\str}{{\rm str}}
\newcommand{\Ind}{{\rm ind}}
\newcommand{\triv}{{\rm triv}}
\newcommand{\Z}{{\mathbb{Z}}}
\newcommand{\bD}{{\bf D}}
\newcommand{\bF}{{\bf F}}
\newcommand{\tX}{{\tt X}}
\newcommand{\Cliff}{{\rm Cliff}}
\newcommand{\tY}{{\tt Y}}
\newcommand{\tZ}{{\tt Z}}
\newcommand{\tV}{{\tt V}}
\newcommand{\tR}{{\tt R}}
\newcommand{\Fam}{{\rm Fam}}
\newcommand{\Cusp}{{\rm Cusp}}
\newcommand{\bT}{{\bf T}}
\newcommand{\bK}{{\bf K}}
\newcommand{\K}{{\mathbb{K}}}
\newcommand{\tH}{{\tt H}}
\newcommand{\bS}{{\bf S}}
\newcommand{\bB}{{\bf B}}
\newcommand{\tW}{{\tt W}}
\newcommand{\tF}{{\tt F}}
\newcommand{\bA}{{\bf A}}
\newcommand{\bL}{{\bf L}}
 \newcommand{\bom}{{\bf \Omega}}
\newcommand{\ch}{{\bf ch}}
\newcommand{\ve}{{\varepsilon}}
\newcommand{\C}{{\mathbb{C}}}
\newcommand{\gen}{{\rm gen}}

\newcommand{\bP}{{\bf P}}
\newcommand{\Naaa}{{\bf N}}
\newcommand{\image}{{\rm image}}
\newcommand{\gaaa}{{\frak g}}
\newcommand{\zaaa}{{\frak z}}
\newcommand{\saaa}{{\frak s}}
\newcommand{\laaa}{{\frak l}}
\newcommand{\stimes}{{\times\hspace{-1mm}\bf |}}
\newcommand{\ausg}{{\rm end}}
\newcommand{\bff}{{\bf f}}
\newcommand{\maaa}{{\frak m}}
\newcommand{\aaaa}{{\frak a}}
\newcommand{\naaa}{{\frak n}}
\newcommand{\brr}{{\bf r}}
\newcommand{\res}{{\rm res}}
\newcommand{\Aut}{{\rm Aut}}
\newcommand{\Pol}{{\rm Pol}}
\newcommand{\Tr}{{\rm Tr}}
\newcommand{\cT}{{\cal T}}
\newcommand{\dom}{{\rm dom}}
\newcommand{\db}{{\bar{\partial}}}
\newcommand{\g}{{\gaaa}}
\newcommand{\cZ}{{\cal Z}}
\newcommand{\cH}{{\cal H}}
\newcommand{\cM}{{\cal M}}
\newcommand{\interi}{{\rm int}}
\newcommand{\singsupp}{{\rm singsupp}}
\newcommand{\cE}{{\cal E}}
\newcommand{\ccR}{{\cal R}}
\newcommand{\cV}{{\cal V}}
\newcommand{\cY}{{\cal Y}}
\newcommand{\cW}{{\cal W}}
\newcommand{\cI}{{\cal I}}
\newcommand{\cC}{{\cal C}}
\newcommand{\cK}{{\cal K}}
\newcommand{\cA}{{\cal A}}
\newcommand{\cEp}{{{\cal E}^\prime}}
\newcommand{\cU}{{\cal U}}
\newcommand{\Hom}{{\mbox{\rm Hom}}}
\newcommand{\vol}{{\rm vol}}
\newcommand{\cO}{{\cal O}}
\newcommand{\End}{{\mbox{\rm End}}}
\newcommand{\Ext}{{\mbox{\rm Ext}}}
\newcommand{\rk}{{\rm rank}}
\newcommand{\im}{{\mbox{\rm im}}}
\newcommand{\sign}{{\rm sign}}
\newcommand{\spann}{{\mbox{\rm span}}}
\newcommand{\symm}{{\mbox{\rm symm}}}
\newcommand{\cF}{{\cal F}}
\newcommand{\cD}{{\cal D}}
\newcommand{\Ree}{{\rm Re }}
\newcommand{\Res}{{\mbox{\rm Res}}}
\newcommand{\Imm}{{\rm Im}}
\newcommand{\inter}{{\rm int}}
\newcommand{\clo}{{\rm clo}}
\newcommand{\tg}{{\rm tg}}
\newcommand{\ee}{{\rm e}}
\newcommand{\Li}{{\rm Li}}
\newcommand{\cN}{{\cal N}}
 \newcommand{\conv}{{\rm conv}}
\newcommand{\op}{{\mbox{\rm Op}}}
\newcommand{\tr}{{\mbox{\rm tr}}}
\newcommand{\cs}{{c_\sigma}}
\newcommand{\ctg}{{\rm ctg}}
\newcommand{\degg}{{\mbox{\rm deg}}}
\newcommand{\Ad}{{\mbox{\rm Ad}}}
\newcommand{\ad}{{\mbox{\rm ad}}}
\newcommand{\codim}{{\rm codim}}
\newcommand{\Gr}{{\mathrm{Gr}}}
\newcommand{\coker}{{\rm coker}}
\newcommand{\id}{{\mbox{\rm id}}}
\newcommand{\ord}{{\rm ord}}
\newcommand{\nat}{{\Bbb  N}}
\newcommand{\supp}{{\rm supp}}
\newcommand{\sing}{{\mbox{\rm sing}}}
\newcommand{\spec}{{\mbox{\rm spec}}}
\newcommand{\Ann}{{\mbox{\rm Ann}}}
\newcommand{\aca}{{\aaaa_\C^\ast}}
\newcommand{\acag}{{\aaaa_{\C,good}^\ast}}
\newcommand{\acage}{{\aaaa_{\C,good}^{\ast,extended}}}
\newcommand{\tck}{{\tilde{\ck}}}
\newcommand{\tnk}{{\tilde{\ck}_0}}
\newcommand{\ceep}{{{\cal E}(E)^\prime}}
 \newcommand{\ncE}{{{}^\naaa\cE}}
 \newcommand{\Or}{{\rm Or }}
\newcommand{\Diff}{{\cal D}iff}
\newcommand{\cB}{{\cal B}}
\newcommand{\hc}{{{\cal HC}(\gaaa,K)}}
\newcommand{\hcma}{{{\cal HC}(\maaa_P\oplus\aaaa_P,K_P)}}
\def\imath{{\rm i}}
\newcommand{\vsl}{{V_{\sigma_\lambda}}}
\newcommand{\czg}{{\cZ(\gaaa)}}
\newcommand{\csl}{{\chi_{\sigma,\lambda}}}
\newcommand{\cR}{{R}}
\def\hB{\hspace*{\fill}$\Box$ \newline\noindent}
\newcommand{\varho}{\varrho}
\newcommand{\ind}{{\rm index}}
\newcommand{\Indu}{{\rm Ind}}
\newcommand{\Fin}{{\mbox{\rm Fin}}}
\newcommand{\cS}{{S}}
\newcommand{\orig}{{\cal O}}
\def\hB{\hspace*{\fill}$\Box$ \\[0.5cm]\noindent}
\newcommand{\cL}{{\cal L}}
 \newcommand{\cG}{{\cal G}}
\newcommand{\npci}{{\naaa_P\hspace{-1.5mm}-\hspace{-2mm}\mbox{\rm coinv}}}
\newcommand{\pki}{{(\paaa,K_P)\hspace{-1.5mm}-\hspace{-2mm}\mbox{\rm inv}}}
\newcommand{\mki}{{(\maaa_P\oplus \aaaa_P, K_P)\hspace{-1.5mm}-\hspace{-2mm}\mbox{\rm inv}}}
\newcommand{\Mat}{{\rm Mat}}
\newcommand{\npi}{{\naaa_P\hspace{-1.5mm}-\hspace{-2mm}\mbox{\rm inv}}}
\newcommand{\ngp}{{N_\Gamma(\pi)}}
\newcommand{\gbg}{{\Gamma\backslash G}}
\newcommand{\gkm}{{ Mod(\gaaa,K) }}
\newcommand{\ggkm}{{  (\gaaa,K) }}
\newcommand{\pkm}{{ Mod(\paaa,K_P)}}
\newcommand{\ppkm}{{  (\paaa,K_P)}}
\newcommand{\makm}{{Mod(\maaa_P\oplus\aaaa_P,K_P)}}
\newcommand{\mmakm}{{ (\maaa_P\oplus\aaaa_P,K_P)}}
\newcommand{\cP}{{\cal P}}
\newcommand{\gm}{{Mod(G)}}
\newcommand{\gk}{{\Gamma_K}}
\newcommand{\La}{{\cal L}}
\newcommand{\cug}{{\cU(\gaaa)}}
\newcommand{\cuk}{{\cU(\kaaa)}}
\newcommand{\dc}{{C^{-\infty}_c(G) }}
\newcommand{\gdk}{{\gaaa/\kaaa}}
\newcommand{\dgkm}{{ D^+(\gaaa,K)-\mbox{\rm mod}}}
\newcommand{\dgm}{{D^+G-\mbox{\rm mod}}}
\newcommand{\vect}{{\C-\mbox{\rm vect}}}
 \newcommand{\cig}{{C^{ \infty}(G)_{K} }}
\newcommand{\gami}{{\Gamma\hspace{-1.5mm}-\hspace{-2mm}\mbox{\rm inv}}}
\newcommand{\cQ}{{\cal Q}}
\newcommand{\mmap}{{Mod(M_PA_P)}}
\newcommand{\bbbz}{{\bf Z}}
 \newcommand{\cX}{{\cal X}}
\newcommand{\bH}{{\bf H}}
\newcommand{\pr}{{\rm pr}}
\newcommand{\bX}{{\bf X}}
\newcommand{\bY}{{\bf Y}}
\newcommand{\bZ}{{\bf Z}}
\newcommand{\bV}{{\bf V}}

\newtheorem{prop}{Proposition}[section]
\newtheorem{lem}[prop]{Lemma}
\newtheorem{ddd}[prop]{Definition}
\newtheorem{theorem}[prop]{Theorem}
\newtheorem{kor}[prop]{Corollary}
\newtheorem{ass}[prop]{Assumption}
\newtheorem{con}[prop]{Conjecture}
\newtheorem{prob}[prop]{Problem}
\newtheorem{fact}[prop]{Fact}

\begin{abstract}
Let $G$ be countable group and $M$ be a proper cocompact
even-dimensional $G$-manifold with orbifold quotient $\bar M$. Let $D$ be a
$G$-invariant Dirac operator on $M$. It induces an equivariant $K$-homology
class $[D]\in K^G_0(M)$ and an orbifold Dirac operator $\bar D$ on $\bar M$.
Composing the assembly map $K^G_0(M)\rightarrow K_0(C^*(G))$
with the homomorphism $K_0(C^*(G))\rightarrow \Z$ given by the representation
$C^*(G)\rightarrow \C$ of the maximal group $C^*$-algebra induced from the trivial representation of $G$ we
define $\ind([D])\in \Z$. In the second  section of the paper we show that
$\ind(\bar D)=\ind([D])$ and obtain explicit formulas for this integer.
In the third section we review the decomposition of
$K_0^G(M)$ in terms of the contributions of fixed point sets of finite cyclic
subgroups of $G$ obtained by W. L\"uck. In particular, the class $[D]$
decomposes in this way.  In the last section we derive an explicit formula for
the contribution to $[D]$ associated to a finite cyclic subgroup of $G$.
\end{abstract}

\tableofcontents

\section{Introduction}

Let $G$ be countable group and $M$ be a proper cocompact
even-dimensional $G$-manifold with orbifold quotient $\bar M$. 
In the literature, orbifolds which can be represented as a global quotient of a smooth manifold by a proper action of a discrete group are often called good orbifolds.

Let $D$ be a
$G$-invariant Dirac operator on $M$ acting on sections of a $G$-equivariant $\Z/2\Z$-graded Dirac bundle $F\to M$. It induces an equivariant $K$-homology
class $[D]\in K^G_0(M)$ and an orbifold Dirac operator $\bar D$ on $\bar M$
with index $\ind(\bar D)\in \Z$.
In the following we briefly describe these objects.

We can identify $\bar D$ with the restriction of $D$ to
the subspace of $G$-invariant sections $C^\infty(M,F)^G$.
The operator $\bar D$ is an example of an elliptic operator on an orbifold. Index theory for elliptic operators on orbifolds has been started with \cite{MR641150} (see also \cite{MR527023}, \cite{MR0474432} for special cases, and \cite{MR1164622}, \cite{MR1163555}, \cite{MR1127139} for alternative approaches). In particular, we have $\dim\:\ker(\bar D)<\infty$, and we can define $$\ind(\bar D):=\dim \:\ker(\bar D^+)-\dim\:\ker(\bar D^-)\ .$$ 

In the present paper we use the analytic definition of equivariant $K$-homology using equivariant $KK$-theory
$$K^G(M):=KK^G(C_0(M),\C)\ .$$
The class $[D]\in KK^G(C_0(M),\C)$ is represented by the Kasparov module
$(\cE,\cF)$ with $\cE:=L^2(M,F)$ and $\cF:=D(D^2+1)^{-1/2}$ (see Subsection \ref{qq12} for more details).


Let $C^*(G)$ denote the unreduced group $C^*$-algebra of $G$. In general, the theory of the present paper would not work with the reduced group $C^*$-algebra $C_r^*(G)$. The key point is that finite-dimensional unitary representations of $G$ extend to representations of $C^*(G)$, but not to $C_r^*(G)$ in general.

We now consider the  assembly map
$$ass:K^G_0(M)\rightarrow K_0(C^*(G))\ .$$
We use an analytic description of the assembly map which is part of Definition \ref{aassww}, and we
refer to \cite{MR2193334}, \cite{MR1659969} and \cite{MR2101228}
for modern treatements of assembly maps in general.

Composing the assembly map  
with the homomorphism $I_1:K_0(C^*(G))\rightarrow K_0(\C)\cong \Z$ given by the representation
$1:C^*(G)\rightarrow \C$ induced from the trivial representation of $G$ we
define 
$$\ind([D]):=I_1\circ ass([D])\in \Z\ .$$

As a special case of the first main result Theorem \ref{orbi}
we get  the equality   
\begin{equation}\label{gghh22}
\ind(\bar D)=\ind([D])\ .
\end{equation}
Theorem \ref{orbi} deals with the slightly more general case where the trivial representation
$triv$ of $G$ is replaced by an arbitrary  finite-dimensional unitary representation of $G$.
We think, that equation (\ref{gghh22})  was  known to specialists, at least as a folklore fact. 

The next result of the present paper is a nice local formula for $\ind([D])$.
The main feature of local index theory is that one can calculate the index of a Dirac operator on a 
closed smooth manifold in terms of an integral of a local index form. A standard reference for local index theory is the book \cite{MR1215720}.  Local index theory generalizes to Dirac operators on orbifolds.
The index formulas in \cite{MR641150} and  \cite{MR1164622} express the index of the Dirac operator on the orbifold as a sum of integrals of local index forms over the various strata. 
In the case of a good orbifold $G\backslash M$ the strata correspond to the fixed point manifolds $M^g$ of the elements $g\in G$. There are various ways to organize these contributions. For the purpose of the present paper we need a formula which expresses the index as a sum of contributions associated to the conjugacy classes of finite cyclic subgroups of $G$.  We will state this formula in Corollary \ref{ccssdd} (we refrain from giving a detailed statement here since this would require the introduction of too much of  notation). In principle one could deduce the formula given in Corollary \ref{ccssdd}  by reorganising the previous results \cite{MR641150} and  \cite{MR1164622}.
But we found it simpler to prove the formula directly using the heat equation approach to local index theory
and the local calculations from equivariant index theory  \cite{MR1215720}.

The proper cocompact $G$-manifold $M$ can be given the structure of a finite $G$-CW-complex.
The equivariant $K$-homology of proper $G$-CW-complexes has been studied intensively in connection
with the Baum-Connes conjecture. Rationally,  $K^G(M)$ decomposes as a sum of contributions of conjugacy classes $(C)$ of finite cyclic subgroups $C\subset G$ (see (\ref{sur}) for a detailed statement). This decomposition is a consequence of a result of \cite{MR1914619} which is finer since it only requires to invert the primes dividing the orders of the finite subgroups of $G$. We thus can write $[D]$ as a sum of contributions $[D](C)$ where $(C)$ runs over the set of conjugacy classes of finite cyclic subgroups of $G$.
Our last result Theorem \ref{hhjj} is the calculation of $[D](C)$. In the proof we use the index formula Corollary \ref{ccssdd} as follows.  By a result of \cite{MR1838997} the equivariant $K$-theory $K_G^0(M)$ has a description in terms of finite-dimensional $G$-equivariant vector bundles $E\to M$. We first derive a cohomological  index formula Theorem \ref{inform} for the pairing of a $K$-homology class coming from a finite cyclic subgroup $C\subset G$ with the class $[E]\in K^0_G(M)$. In the proof we use the relation (\ref{gghh22}).

We then observe that the pairing of $[D]$ with $[E]$ is the index of the twisted operator
$[D_E]$ which can be written as a sum of contributions of conjugacy classes of finite subgroups by  \ref{ccssdd}.
We obtain $[D](C)$ be a comparison of the formulas  in Theorem \ref{inform} and Corollary \ref{ccssdd}
and variation of $E$.

{\em Acknowledgement: The first version of this paper was written in spring 2001. I want to thank W. L{\"u}ck for his motivating interest in this work, and Th. Schick for pointing out a small mistake \footnote{The factor $\frac{1}{\ord(g)}$ in (\ref{rree33}) was missing.} in the previous version.}

\section{Assembly and orbifold index}

\subsection{The equivariant $K$-homology class of an invariant Dirac operator}\label{qq12}

Let $G$ be a countable discrete group. Let $M$ be a smooth proper cocompact
$G$-manifold, i.e.  a $G$-manifold such that the stabilizer $G_x$ is finite
for all $x\in M$, and $G\backslash M$ is compact. We further assume that $M$
is equipped with a complete $G$-invariant Riemannian metric $g^M$ and a
$G$-homogeneous Dirac bundle $(F,\nabla^F,\circ,(.,.)_F)$. Here
$\circ:TM\otimes F\rightarrow F$ is the Clifford multiplication, $\nabla^F$ is
a Clifford connection,  $(.,.)_F$ is the hermitian scalar product, and these
structures satisfy the usual compatibiliy conditions (see
\cite{MR1215720}, Ch.3)  and are, in addition, $G$-invariant.

For simplicity we assume that $\dim(M)$ is even and that the Dirac bundle is
$\Z/2\Z$-graded. In fact, the odd-dimensional case can easily be reduced to the
even dimensional case by taking the product with $S^1$.

We use equivariant $KK$-theory in order to define equivariant $K$-homology.
Thus let $KK^G$ be the equivariant $KK$-theory introduced in \cite{MR918241} (see also  \cite{MR1656031}).
Let $C_0(M)$ be the $G$-$C^*$algebra of continuous functions on $M$ vanishing
at infinity. Then by definition $K^G_0(M)=KK^G(C_0(M),\C)$.
The Dirac operator $D$ associated to the invariant Dirac bundle $F$ induces a
class $[D]\in K^G_0(M)$ as follows. We form the $\Z/2\Z$-graded $G$-Hilbert space
$\cE:=L^2(M,F)$. Then $C_0(M)$ acts on $\cE$ by multiplication. Furthermore,
we consider the bounded $G$-invariant operator  $\cF:=D(D^2+1)^{-1/2}$
which is defined by applying the function calculus to the unique (see
\cite{MR0369890}) selfadjoint extension of $D$. Then $[D]$ is represented
by the Kasparov module $(\cE,\cF)$.  

\subsection{Descent and index}

Let $C^*(G)$ denote the (non-reduced) group $C^*$-algebra of $G$. It has the
universal property, that any unitary representation of $G$ extends to
representation of $C^*(G)$. In particular, if $\rho:G\rightarrow U(V_\rho)$
is an unitary representation of $G$ on  a finite-dimensional Hilbert space $V_\rho$, then there is an extension
$\rho:C^*(G)\rightarrow \End(V_\rho)$. On the level of $K$-theory it induces a
homomorphism (using Morita invariance and $K_0(\C)\cong \Z$)
$I_\rho:K_0(C^*(G))\rightarrow K_0(\End(V_\rho))\cong \Z$. 
In particular, if $\rho=1$ is the trivial representation, then we also write
$I:=I_1$. Note that $I_\rho$ can be written as a Kapsarov
product $\otimes_{C^*(G)} [\rho]$, where $[\rho]\in KK(C^*(G),\End(V(\rho)))$
is represented by the Kasparov module $(V_\rho,0)$.

Let $C^*(G,C_0(M))$ be the (non-reduced) cross product of $G$ with $C_0(M)$.
Then there is the descent homomorphism
 $j^G:K^G_0(M)\cong KK^G(C_0(M),\C)\rightarrow KK(C^*(G,C_0(M)),C^*(G))$
introduced in \cite{MR918241}, 3.11.
Following \cite{MR1711324} we choose any cut-off function $\chi\in
C^\infty_c(M)$ with values in $[0,1]$ such that $\sum_{g\in G} g^*\chi^2\equiv
1$. Then we define the projection $P\in C^*(G,C_0(M))$ by $P(g)=(g^{-1})^*\chi
\chi$. Let $[P]\in K_0(C^*(G,C_0(M))\cong KK(\C, C^*(G,C_0(M)))$ be the class
induced by $P$, which is independent of the choice of $\chi$. 

\begin{ddd}\label{aassww}
We define 
$\ind_\rho :K^G_0(M)\rightarrow \Z$
to be the composition
$$K^G_0(M)\stackrel{j^G}{\rightarrow}KK(C^*(G,C_0(M)),C^*(G))\stackrel{[P]\otimes_{C^*(G,C_0(M))}}{\longrightarrow} KK(\C, C^*(G,C_0(M)))
\stackrel{I_\rho}{\rightarrow} \Z\ .$$
\end{ddd}

In particular, we set $\ind:=\ind_1$.

\subsection{Index and Orbifold index}

The quotient $\bar M:=G\backslash M$ is a smooth compact orbifold carrying an
orbifold Dirac bundle $\bar F:=G\backslash F$ with associated orbifold Dirac
operator $\bar D$. In our case the space of smooth sections  $C^\infty(\bar
M,\bar F)$ can be identified with the $G$-invariant sections
$C^\infty(M,F)^G$. Then $\bar D$ coincides with the restriction of $D$ to this
subspace. It is well-known that $\dim(\ker \bar D)<\infty$ so that we can
define the index $\ind(\bar D):=\dim_s\ker(\bar D)\in \Z$, where the
subscript "${}_s$" indicates hat we take the super dimension.

If $\rho:G\rightarrow U(V_\rho)$ is a finite-dimensional unitary representation
of $G$, then we define the orbifold bundle $\bar V(\rho):=G\backslash M\times
V_\rho$ and let $\bar D_\rho$ be the twisted operator associated to $\bar
F\otimes \bar V(\rho)$.  The space $C^\infty(\bar M,\bar
F\otimes \bar V(\rho))$ can be identified with $(C^\infty(M,F)\otimes
V_\rho)^G$ such that $\bar D_\rho$ is the restriction of $D\otimes 1$ to this
subspace. Still we can define
$\ind(\bar D_\rho)$.

\begin{theorem}\label{orbi}
$\ind(\bar D_\rho)=\ind_\rho([D])$ 
\end{theorem}
\proof
We first apply $j^G$ to the Kasparov module $(L^2(M,F),\cF)$ representing
$[D]$. According to \cite{MR918241}, 3.11., $j^G([D])$ is represented by
$(C^*(G,L^2(M,F)),\tilde\cF)$, where $C^*(G,L^2(M,F))$
is a $C^*(G)$-right-module admitting a left action by $C^*(G,C_0(M))$.
It is a closure of the space of finitely supported functions $f:G:\rightarrow
L^2(M,F)$. The operator $\tilde \cF$ is given by $(\tilde\cF f)(g)=(\cF f)(g)$.
The $C^*(G)$-valued scalar product is given by $\langle
f_1,f_2\rangle(g)=\sum_{h\in G}\langle f_1(h), f_2(hg)\rangle $. 
Furthermore, the left action of $C^*(G,C_0(M))$ is given by $(\phi
f)(g)=\sum_{h\in G} \phi(h) (hf)(g)$.

Using associativity of the Kasparov product we can compute $\ind_\rho$ by
first applying $\otimes_{C^*(G)} [\rho]$ and then $[P]\otimes_{C^*(G,C_0(M))}$.
Using that $C^*(G,L^2(M,F))\otimes _{C^*(G)}
V_\rho\cong L^2(M,F)\otimes V_\rho$ by
$f\otimes v\mapsto  \sum_{g\in G}  f(g) \rho(g)v$
 we conclude that
$j^G([D])\otimes_{C^*(G)} [\rho]$ is represented by the Kasparov module
$(L^2(M,F)\otimes V_\rho,\hat \cF)$, where $\hat \cF=\cF\otimes
\id_{V_\rho}$. The left-action of $C^*(G,C_0(M))$ is given by
$(\phi f)=\sum_{h\in G} \phi(h) (h\otimes \rho(h) )f$.

Finally we compute $[P]\otimes_{C^*(G,C_0(M))}\left(j^G([D])\otimes_{C^*(G)}
[\rho]\right)$. We represent $[P]$ by the Kasparov module
$(PC^*(G,C_0(M)),0)$. We must understand $PC^*(G,C_0(M))
\otimes_{C_0(M)} (L^2(M,F)\otimes V_\rho)$.

There is a natural unitary inclusion $L:L^2(\bar M,\bar F\otimes \bar
V(\rho))\hookrightarrow L^2(M,F)\otimes V_\rho$. If $f\in L^2(\bar M,\bar F\otimes \bar
V(\rho))$ is considered as an element $\hat f$ of $(L^2_{loc}(M,F)\times
V_\rho)^G$ in the natural way, then $L(f):=\chi \hat f$.
The projection $LL^*$ onto the range of $L$ is given by
$$LL^*(f)=\sum_{g\in G} (g^{-1})^*\chi  g f\ .$$
It now follows from the definition of $P$ that
\begin{eqnarray*}
PC^*(G,C_0(M))
\otimes_{C^*(G,C_0(M))}(L^2(M,F)\otimes V_\rho)&=& P (L^2(M,F)\otimes V_\rho)\\
&\stackrel{L^*}{\cong}& L^2(\bar M,\bar F\otimes \bar
V(\rho))
\end{eqnarray*}
The operator $\bar D$ has a natural selfadjoint extension (also denoted by
$\bar D$ such that we can form $\bar \cF:=\bar D(1+\bar D^2)^{-1/2}$.
We claim that $[P]\otimes_{C^*(G,C_0(M))}\left(j^G([D])\otimes_{C^*(G)}
[\rho]\right)$ is represented by the Kapsarov module
$(L^2(\bar M,\bar F\otimes \bar
V(\rho)),\bar \cF)$. The assertion of the Theorem immediately follows from the
claim. In order to show the claim we employ the characterization of the Kasparov product in terms of
connections (see \cite{MR918241}, 2.10). In our situation we have only to
show that $\bar\cF$ is a $\hat \cF$-connection.

For Hilbert-$C^*$-modules $X,Y$ over some $C^*$-algebra $A$ let $L(X,Y)$ and $K(X,Y)$ denote the spaces of bounded and compact adjoinable $A$-linear operators (see \cite{MR1656031} for definitions).
For $\xi\in PC^*(G,C_0(M))$ we define $\theta_\xi\in
L(L^2(M,F)\otimes V_\rho, PL^2(M,F)\otimes
V_\rho )$ by $\theta_\xi(f)=\xi f$.
Since $\cF$ and $\bar\cF$ are selfadjoint   we only must show that
$\theta_\xi\circ \hat \cF- (L\bar \cF L^*) \circ \theta_\xi\in 
K(L^2(M,F)\otimes V_\rho, PL^2(M,F)\otimes V_\rho )$.  
We have $\xi \hat \cF- (L\bar \cF L^*) \xi
= [\xi,\hat \cF]+(\hat \cF -L\bar \cF L^*)  P\xi$.
Since $[\xi,\hat \cF]$ is compact it suffices to show that $(\hat \cF -(L\bar
\cF L^*)  P$ is compact. 
We consider $\tilde D:=(1-P)D(1-P)+L\bar D L^*$. 
Then we have $\tilde D=D+Q$, where $Q$ is a  zero order non-local operator.
Let $\tilde\cF:=\tilde \cD(1+\tilde \cD^2)^{-1/2}$. Then
$(\hat \cF -L\bar \cF L^*) P=(\hat \cF-\tilde \cF)P$.
Let $\tilde \chi\in C_c^\infty(M)$ be such that $\chi\tilde\chi=\chi$.
Then we have $(\hat \cF-\tilde \cF)P=(\hat \cF-\tilde \cF)\tilde\chi P$.
Therefore it suffices to show that $(\hat \cF-\tilde \cF)\tilde\chi$ is
compact. This can be done using the integral representations for $\hat \cF$
and $\tilde \cF$ as in \cite{MR1348799}.
\hB

\subsection{The local index theorem}

In this the present subsection we derive a local index theorem which is
a formula for $\ind_\rho([D])$ in terms of integrals of characteristic forms
over the various singular strata of $\bar M$.  

Let $W\in C^\infty(M\times M, F\boxtimes F^*)^G$ be a an invariant section which
satisfies an estimate 
\begin{equation}\label{est}
|W(x,y)|\le C\exp(-c\dist(x,y)^2)\end{equation}
 for some $c>0$,
$C<\infty$. Since $\bar M$ is compact the manifold $M$ has bounded geometry,
and in particular,  it has at most exponential volume growth. Therefore, $W$
defines an integral operator $\bar W$  on $L^2(\bar M,\bar F\otimes \bar
V_\rho)$ by $$\bar Wf(x):=\int_M (W(x,y)\otimes \id_{V_\rho})  f(y) dy\ .$$
This operator is in fact of trace class. We claim that
\begin{equation}\label{tre}\Tr \:\bar W= \int_{\bar M}\sum_{g\in G} \tr (W(x,g
x) g_x)  dx \:\tr \rho(g)\ ,\end{equation}
where $g_x$ denotes the linear map $g_x:F_x\rightarrow
F_{gx}$.
In order to see the claim note that $\Tr\: \bar W= \Tr \:LW\bar L^*$, and
$R:=LW  L^*$ is the integral operator on $L^2(M,F)\otimes V_\rho$ given by
the integral kernel  $R (x,y)=\sum_{g\in G} \chi(x) W(x,g y)  g_y  \chi (y)
\otimes \rho(g) $.

Again, since $M$ and $F$ have bounded geometry the heat kernel
$W_t$, $t>0$, i.e. the integral kernel of $\exp(-tD^2)$, satisfies the
Gaussian estimate (\ref{est}). Moreover, $\bar W_t$ is precisely $\exp(-t\bar
D^2)$. By the McKean-Singer formula we have
$$\ind(\bar D_\rho)=\Tr_s \bar W_t$$ 
for any $t>0$,  where $\Tr_s$ is the super trace. We obtain the local index
formula by evaluating $\lim_{t\to 0} \Tr_s \bar W_t$.

If $g\in G$, then let $M^g$ denote the fixed point submanifold of $g$.
If $M^g\not=\emptyset$, then $g$ is of finite order.
Furthermore, let  $Z_G(g)$ denote the centralizer of $g$ in $G$.
Then $Z_G(g)\backslash M^g$ is compact.
For $g\in G$ let  $(g)\in C(G)$
denote the conjugacy class of $g$, where $C(G)$ denotes the set of conjugacy
classes.
By $\cF(G)$ we denote the set of elements of finite order, and by $\cF C(G)$
we denote the set of conjugacy classes of $G$ of finite order.

The formula (\ref{tre}) can we rewritten as follows.
\begin{eqnarray*}
\Tr_s \bar W&=&\int_{\bar M} \sum_{g\in G} \tr_s (W(x,g x)
g_x)  dx \:\tr \rho(g)\\
&=&\sum_{(g)\in C(G)} \int_{ G\backslash M}\sum_{h\in
Z_G(g)\backslash G}  \tr_s (W(x,hgh^{-1} x)
(hgh^{-1})_x)  dx \: \tr \rho(hgh^{-1})\\
&=&\sum_{(g)\in C(G)}  \int_{ Z_G(g)\backslash M}  \tr_s (W(x,g x)
g_x)  dx \: \tr \rho(g)\ .
\end{eqnarray*}
If $W=W_t$ is the heat kernel, then due to the usual gaussian estimates the integral $\int_{ Z_G(g)\backslash M} 
\tr_s (W(x,g x) g_x)  dx$ localizes at $Z_G(g)\backslash M^g$ as $t\to 0$.
There is a $Z_G(g)$-invariant density
$U(g)\in C^\infty(M^g, |\Lambda^{max}|T^*M^g)^{Z_G(g)}$ which is locally
determined by the Riemannian structure $g^M$ and the Dirac bundle $F$
such that 
$$\lim_{t\to\infty} \int_{ Z_G(g)\backslash M} 
\tr_s (W_t(x,g x) g_x)  dx=\frac{1}{\ord(g)}\int_{Z_G(g)\backslash M^g} U(g)\ .$$
An explicit formula for $U(g)$ is given in \cite{MR1215720}, Ch.
6.4, and it will be recalled below.
We conclude that
$$\ind_\rho([D])=\sum_{(g)\in C\cF (G)} \frac{1}{\ord(g)} \int_{Z_G(g)\backslash M^g} U(g) \tr
\rho(g)\ .$$

The fixed point manifold $M^g$ is a totally geodesic Riemannian submanifold
of $M$ with induced metric $g^{M^g}$. Let $R^{M^g}$ denote its  curvature
tensor. We define the form $\hat \bA(M^g)\in \Omega(M^g,\Or(M^g))$ by
$$\hat\bA(M^g)={\det}^{1/2}\left(\frac{R^{M^g}/4\pi\imath}{\sinh(R^{M^g}/4\pi\imath)}\right)\
,$$ where $\Or(M^g)$ denote the orientation bundle (the orientation bundle
occurs since we must choose an orientation in order to define $\det^{1/2})$.

Furthermore, we define the $G$-equivariant bundle $F/S:=\End_{\Cliff(TM)}(F)$.
It comes with a natural connection $\nabla^{F/S}$. By $R^{F/S}$ we denote its
curvature. Following \cite{MR1215720}, 6.13, we define the form
$\ch(g,F/S)\in \Omega(M^g, \Lambda^{max} N\otimes\Or(M))$ by
$$\ch(g,F/S)=\frac{2^{\codim_M(M^g)}}{\sqrt{\det(1-g^N)}} \str
(\sigma_{\codim_M(M^g)}(g^F) \exp (-R_0^{F/S}/2\pi\imath))\ .$$ 
Here $g^N$ is the restriction of $g$ to the normal bundle $N$ of $M^g$. Note
that $\det(1-g^N)>0$ so that $\sqrt{\det(1-g^N)}$ is well-defined.
Furthermore $g^F$ is the action of $g$ on the fibre of $F_{|M^g}$.
Since $g^F$ commutes with $\Cliff(TM^g)$ it corresponds to an element of
$\Cliff(N)\otimes \End_{\Cliff(M)}(F)$. 
$\sigma_{\codim_M(M^g)}:\Cliff(N)\rightarrow \Lambda^{max} N$ is the symbol map
so that $  \sigma_{\codim_M(M^g)}g^F\in \End_{\Cliff(M)}(F)\otimes
\Lambda^{max} N$. Furthermore, the restriction $R_0^{F/S}$ of the curvature
$R^{F/S}$ to $M^g$ is a section of
$\Omega(M^g,\End_{\Cliff(M)}(F)_{|M^g})$. The super trace
$\str:\End_{\Cliff(M)}(F)\rightarrow \C\otimes \Or(M)$ is defined by
$\str(W)=\tr_s(\Gamma W)$, where $\Gamma=\imath^{n/2} \vol_M$ is the chirality
operator defined using the orientation of $M$.

Let
$T_N: \Lambda^{max}N\rightarrow \C\otimes \Or(N)$ be the normal Beresin
integral, where $\Or(N)$ is the bundle of normal orientations. Then
we have
$$\label{ug}U(g):=[T_N(\frac{\hat
\bA(M^g)\ch(g,F/S)}{\det^{1/2}(1-g^N\exp(-R^N/2\pi\imath))})]_{max}\
.$$  Here $R^N$ is the curvature tensor of $N$,
$\frac{1}{\det^{1/2}(1-g^N\exp(-R^N))}\in \Omega(M^g,\Or(M^g))$, and
$[.]_{max}$ takes the part of maximal degree. In order to interpret the
right-hand side as a density on $M^g$ we identify $\Lambda^{max}T^*M^g\otimes
\Or(M^g)^2\otimes\Or(N) \otimes \Or(M)$ with $|\Lambda^{max}|T^*M^g$ in the
canonical way.

\begin{theorem}\label{rree33}
$$\ind_\rho([D])=\sum_{(g)\in C\cF (G)} \frac{\tr
\rho(g)}{\ord(g)}   \int_{Z_G(g)\backslash M^g} [T_N(\frac{\hat
\bA(M^g)\ch(g,F/S)}{\det^{1/2}(1-g^N\exp(-R^N/2\pi\imath))})]_{max} $$
\end{theorem}

\subsection{Cyclic subgroups}

We now reformulate the local index theorem in terms of contributions of
conjugacy classes of cyclic subgroups.
Let $\cF Cyc(G)$ denote the set  of finite cyclic
subgroups. If $C\in \cF Cyc(G)$, then let $\gen(C)$
denote the set of its generators. The normalizer $N_G(C)$ and the Weyl group
$W_G(C):=N_G(C)/Z_G(C)$ acts on $\gen(C)$.
There is a natural map $p:\cF(G)\rightarrow \cF Cyc(G)$, $g\mapsto <g>$ which
factors over conjugacy classes $\bar p:C\cF(G)\rightarrow C\cF Cyc(G)$.
If $(C)\in C\cF Cyc(G)$, then $\bar p^{-1}(C)$ can be identified with 
 $W_G(C)\backslash \gen(C)$.

Note that $M^g=M^{<g>}$, i.e. it only depends on the cyclic subgroup generated
by $g$. Similarly, $Z_G(g)=Z_G(<g>)$.  
So we obtain
\begin{kor}\label{ccssdd}
$$\ind_\rho([D])=\sum_{(C)\in C \cF Cyc (G)} \frac{1}{|C|} \sum_{g\in
W_G(C)\backslash \gen(C)} \int_{Z_G(C)\backslash M^C} U(g) \tr \rho(g) $$
\end{kor}

\subsection{Cap product and twisting}

We define $K_G^0(M):=KK^G(\C,C_0(M))$. If
$E$ is a $G$-equivariant complex vector bundle, then
let $[E]\in K_G^0(M)$ denote the class  represented by
the Kasparov module $(C_0(M,E),0)$, where we define the $C_0(M)$-valued scalar
product on $C_0(M,E)$ after choosing a $G$-invariant hermitean metric
$(.,.)_E$.

Since $C_0(M)$ is commutative any right $C_0(M)$-module is a left-
$C_0(M)$-module in a natural way. If we apply this to Kasparaov modules we
obtain a map $$a:KK^G(\C,C_0(M))\rightarrow KK^G(C_0(M),C_0(M))\ .$$
\begin{ddd}
The cap-product $K^0_G(M)\otimes K_0^G(M)\rightarrow K_0^G(M)$ is defined by
$$v\cap x:= a(v)\otimes_{C_0(M)} x\ .$$
\end{ddd}
 
If we choose on $(E,(.,.))$ a hermitian connection $\nabla^E$, then we can
form the twisted Dirac bundle $E\otimes F$ with associated Dirac operator
$D_E$. The following fact is well-known. An elementary proof (for trivial
$G$) can be found e.g. in \cite{MR1348799}.
\begin{prop} \label{cca}
$[D_E]=[E]\cap [D]$
\end{prop}

\subsection{A cohomological index formula for twisted operators}

Let $R^E$ denote the curvature of the connection $\nabla^E$. For a finite cyclic subgroup
$C\subset G$ let $R^E_0$ denote the restriction of $R^E$ to $M^C$.
If $g\in\gen(C)$ , then we have 
$$\ch(g,E\otimes F/S)=\ch(g,F/S) \cup \ch(g,E)\ ,$$
where $\ch(g,E)=\tr g^E \exp(-R_0^E/2\pi\imath)$.
Here $g^E$ denotes the action of $g$ on the fibre of $E$.
Thus we can write
$$U_E(g):=[T_N(\frac{\hat
\bA(M^g)\ch(g,F/S)\cup \ch(g,E)}{\det^{1/2}(1-g^N\exp(-R^N/2\pi\imath))})]_{max}\ .$$
We can write 
$U_E(g)=[\hat U(g) \cap \ch(g,E)]_{max}$,
where 
\begin{equation}\label{hatudef}
\hat U(g)=T_N(\frac{\hat
\bA(M^g)\ch(g,F/S)}{\det^{1/2}(1-g^N\exp(-R^N/2\pi\imath))})\ .
\end{equation}

The cohomology $H^*( Z_G(C)\backslash M^C,\C)$ of the orbifold
$Z_G(C)\backslash M^C$ can be computed using the complex of invariant
differential forms  $(\Omega^*(M^C)^{Z_G(C)},d)$. Furthermore, the homology
$H_*(Z_G(C)\backslash M^C,\C)$ can be identified with the dual of the
cohomology, i.e.   $H_*(Z_G(C)\backslash M^C,\C)\cong H^*( Z_G(C)\backslash
M^C,\C)^*$. The closed form $\hat U(g)\in \Omega^*(M^g,\Or)$ now defines a
homology class $[\hat U(g)]\in H_*(Z_G(C)\backslash M^C,\C)$ such that
$[\hat U(g)]([\omega])=\int_{Z_G(C)\backslash M^C}  [\hat U(g)\cap
\omega]_{max}$ for any closed form  $\omega\in \Omega^*(M^C)^{Z_G(C)}$.

Let $[\ch(g,E)]\in H^*(Z_G(C)\backslash M^C,\C)$ denote the cohomology class
represented by the closed form $\ch(g,E)$.

\begin{theorem}\label{oopp}
$$\ind_\rho([E]\cap [D])=\sum_{(C)\in C \cF Cyc (G)}  \frac{1}{|C|}
\sum_{g\in W_G(C)\backslash \gen(C)} \langle [\ch(g,E)], [\hat
U(g)]\rangle   \tr \rho(g) $$ 
\end{theorem}

\section{Chern characters}

\subsection{The cohomological Chern character}

In this Subsection we review the construction of the Chern character given in
\cite{MR1851256}. There the equivariant $K$-theory is introduced
using a classifying space $\bK_G\C$. If $X$ is a proper $G$-CW complex, then
$\bK^0_G(X):=[X,\bK_G\C]_G$, where $[.,.]_G$ denotes the set of
homotopy classes of equivariant maps.

Let $\K_G(X)$ be the Grothendieck group of $G$-equivariant complex vector
bundles. Then  there is a natural homomorphism
$b:\K_G(X)\rightarrow \bK^0_G(X)$,  which is an isomorphism if
$X$ is finite (\cite{MR1851256}, Prop. 1.5).

If $H$ is a finite group, then let $R_\C(H)$ denote the complex representation
ring of $H$ with complex coefficients. The character gives a natural
identification of $R_\C(H)$ with the space of complex-valued class functions
on $H$, i.e. $\C(C(H))$.

Since we want to work with differential forms later on we simplify matters by
working with complex coefficients (the constructions in \cite{MR1851256}
are finer since they work over $\Q$). For any finite subgroup $H\subset G$
the construction \cite{MR1851256}, (5.4),  provides a
homomorphism 
$$\ch^H_X:K_G^0(X)\rightarrow H^*(Z_G(H)\backslash X^H)\otimes \C(C(H))\ .$$
For our purpose it suffices to understand
$\ch^H_X(b(\{E\}))$, where $E$ is a $G$-equivariant complex vector bundle over
$X$, and $\{E\}$ denotes its class in $\K_G(X)$. 
First of all note that $E_{|X^H}$ is a $N_G(H)$-equivariant bundle over
$X^H$. We can further write $E_{|X^H}=\sum_{\phi\in \hat
H}\Hom_H(V_\phi,E_{|X^H})\otimes V_\phi$, where $\Hom_H(V_\phi,E_{|X^H})$ is a
$Z_G(H)$-equivariant bundle over $X^H$. We therefore obtain an element of
$\K_{Z_G(H)}^0(X^H)\otimes R(H)$. We now apply the composition
\begin{eqnarray*}\lefteqn{\K_{Z_G(H)}^0(X^H)\stackrel{\pr^*}{\rightarrow}
\K_{Z_G(H)}^0(EG\times X^H)\stackrel{\cong}{\rightarrow}
\K^0_1(EG\times_{Z_G(H)} X^H) }&&\\&&\stackrel{\ch}{\rightarrow}
H^{*}(EG\times_{Z_G(H)} X^H,\C)\stackrel{(\pr^*){-1}}{\rightarrow}
H^*(Z_G(H)\backslash X^H,\C)\end{eqnarray*}  to the first component, and the
character $R(H)\rightarrow \C (C(H))$ to the second.  The result belongs to
$H^*(Z_G(H)\backslash X^H,\C)\otimes \C(C(H))$ and is $\ch^H_X(b(\{E\}))$.

If $C$ is a finite cyclic subgroup, then let $r:\C(C(C))\rightarrow
\C(\gen(C))$ be the restriction map. 
Note that $W_G(C)$ acts on $\C(\gen(C))$ as well as on $H^*(Z_G(C)\backslash
X^C,\C)$. 
The result \cite{MR1851256}, Lemma 5.6, now asserts that if $X$
is finite, then
\begin{equation}\label{sur1}\prod_{(C)\in C \cF Cyc(G)}(1\otimes r) \ch_X^C:
\bK^0_G(X)_\C\rightarrow \prod_{(C)\in C\cF Cyc(G)}
\left(H^{ev}(Z_G(C)\backslash X^C,\C)\otimes
\C(\gen(C))\right)^{W_G(C)}\end{equation} is an isomorphism.

\subsection{Differential forms}

In the present subsection we give a description of the equivariant Chern
character using differential forms. Let $M$  be a smooth proper $G$-manifold
and $E$ be a $G$-equivariant complex vector bundle over $M$. Then we can find a
$G$-invariant hermitian metric $(.,.)_E$ and a $G$-invariant metric connection
$\nabla^E$. Let $R^E$ denote the curvature of $\nabla^E$. We define the closed
$G$-invariant form $\ch(E)\in \Omega(M)^G$ by $\ch(E):=\tr
\exp(-R^E/2\pi\imath)$. It represents a cohomology class $[\ch(E)]\in H^*(G\backslash
M,\C)$.  Furthermore, we have the class $\ch^{\{1\}}_M(b(\{E\}))$, which is
given by the following composition 
\begin{eqnarray}\lefteqn{\K_{G}^0(M)\stackrel{\pr_1^*}{\rightarrow}
\K_{G}^0(EG\times M)\stackrel{\cong}{\rightarrow} \K_1^0(EG\times_GM
) }&&\nonumber \\&&\stackrel{\ch}{\rightarrow}
H^{*}(EG\times_GM,\C)\stackrel{(\pr_2^*)^{-1}}{\rightarrow} H^*(G\backslash
M,\C) \label{comp}\ .\end{eqnarray}
\begin{lem}\label{notw}
$[\ch(E)]= \ch^{\{1\}}_M(b(\{E\}))$
\end{lem}
\proof
We show that $\ch^{\{1\}}_M(b(\{E\}))$
can be represented by the form  $\ch(E)$.
To do so we employ an approximation $j:\tilde EG\rightarrow EG$,
where $\tilde EG$ is a free $G$-manifold and the $G$-map $j$ is 
$\dim(M)+1$-connected.  This existence of such approximations will be shown
in Subsection \ref{nanana}.
Then we can define $ \ch^{\{1\}}_M(b(\{E\}))$ by
(\ref{comp}) but with $EG$ replaced by $\tilde EG$.
It is now clear that $\pr^*_2 \ch(E)=\ch(G\backslash \pr_1^* E)$.
\hB

Let $C\subset G$ be
a finite cyclic subgroup. Furthermore, let $[\ch(g,E)]\in H^*(Z_G(C)\backslash
M^C,\C)$ denote the cohomology class represented by $\ch(g,E)$. The function
$\gen(C)\ni g\mapsto [\ch(g,E)]$ can naturally be considered as an element 
$[\ch(.,E)] \in H^*(Z_G(C)\backslash M^C,\C)\otimes
\C(\gen(C))$ which is in fact $W_G(C)$-equivariant.
\begin{prop}\label{chern}
$[\ch(.,E)]=(1\otimes r) \ch^C_M(b(\{E\}))$.
\end{prop}
\proof
 First of all note that $R_0^E$ is the curvature of $E_{|M^C}$.
Furthermore, the decomposition $E_{|M^C}=\sum_{\phi\in \hat
C}E(\phi) \otimes V_\phi$ is preserved by $R_0^E$, where
$E(\phi)=\Hom_C(V_\phi,E_{|M^C})$. Let $R^{E(\phi)}$ be the restriction of the
curvature to the subbundle $E(\phi)\otimes V_\phi$. We get for $g\in\gen(C)$
\begin{eqnarray*}(1\otimes r) \ch^C_M(b\{E\})(g)
&\stackrel{def.}{=}&\sum_{\phi\in \hat C}\ch^{\{1\}}_{M^C}(b\{E(\phi)\})\tr\phi(g)
\\&\stackrel{Lemma \ref{notw}}{=}&\sum_{\phi\in \hat C}
 [\ch(E(\phi))] \tr\:\phi(g)\\&=&\sum_{\phi\in \hat C} [\tr
\exp(-R^{E(\phi)}/2\pi\imath)] \tr\phi(g)
\\&=&[\tr g^E \exp(-R_0^E/2\pi\imath)] \\
&=& [\ch(g,E)]\ .\end{eqnarray*}
\hB

\subsection{Smooth approximations of $CW$-complexes}\label{nanana}

The goal of this subsection is to show  that the
approximation $j:\tilde EG\to EG$ used in the proof of Lemma \ref{notw} exists. 
We start with the following general result.
\begin{prop}\label{nasowas}
If $X$ is a countable finite-dimensional $CW$-complex, then there exists a smooth manifold
$M$ and a homotopy equivalence $M\stackrel{\sim}{\to} X$.
\end{prop}
\proof
Let $X$ be a finite-dimensional $CW$-complex.
Following \cite{MR0142128} we call a manifold with boundary $(\bar M,\partial \bar M)$ a tubular
neighbourhood of $X$ if there exists a continuous map $F:\partial \bar M\to X$ such that the underlying topological space of 
$\bar M$  is the mapping cylinder $C(F)=\partial \bar M\times [0,1]\cup_F X$ of $F$, 
the inclusion $\partial \bar M\times [0,1)\hookrightarrow M$ is smooth, and the inclusion
$X\hookrightarrow \bar M$ is smooth on each open cell of $X$.  

Let $$X^0\subseteq X^1\subseteq X^2\subseteq\dots\subseteq X^n=X$$
be the filtration of $X$ by skeletons. We obtain $X^{i+1}$ from $X^i$ by attaching a countable
number of $i+1$-cells.

We will construct a tubular neighbourhood $(\bar M,\partial \bar M)$ 
of $X$ by induction over the dimension with the special property that the tangent bundle $T\bar M\to \bar M$ is trivial. It is clear that the collection of points $X^0$ admits such a tubular neighbourhood of any given dimension. Let us now assume that
there exists an $i$-dimensional  $CW$-complex $Y$ together with a homotopy equivalence $h:Y\stackrel{\sim}{\to} X^i$ and a $m$-dimensional tubular neighbourhood $(W,\partial W)$, $F:\partial W\to Y$ of $Y$, such that $m\ge 2n+3$ and $TW\to W$ is trivial.

Via a homotopy inverse of $h$ the attaching data for $X^i\subseteq X^{i+1}$ yields the attaching data $\tilde \chi_\alpha:S^{i}\to Y$, $\alpha\in J$, for a countable collection $J$ of $i+1$-cells. Let $\tilde Z$ be the result of attaching these cells to $Y$. Then we have a homotopy equivalence $\tilde Z\stackrel{\sim}{\to} X^{i+1}$.

We consider $\R$ with the cell-structure given by its decomposition into unit intervals.
The product
$(W\times \R,\partial W\times \R)$ is a tubular neighbourhood of $Y\times \R$ in a natural way with retraction
$F\times \id_\R:\partial W\times \R\to Y\times \R$.
The projection $Y\times \R\to Y$ is a homotopy equivalence. 

We now fix an inclusion $J\subseteq\Z$ and define attaching maps
$\chi_\alpha:=\tilde\chi_\alpha\times \{\alpha\}:S^i\to Y\times \R$.
Let $\hat Z$ denote the complex obtained by attaching the cells to $Y\times \R$.
Our choice of attaching maps is made such that these $i+1$-cells are attached to the $i$-skeleton of $Y$.
We have a homotopy equivalence $\hat Z\stackrel{\sim}{\to}\tilde Z$.

In order to improve the attaching maps we argue as in the proof of \cite[Theorem II]{MR0142128}.
Since $2\dim(Y)+1=2i+1\le 2n+3\le m=\dim(W)$ we can deform the attaching map $\tilde \chi_\alpha$
slightly so that its image is disjoint from $Y$. 
To do so we adapt the method of the proof of \cite[Theorem 11a]{MR0073980}, and we use the assumption that the open cells of $Y$ are smoothly embedded. 
Using the mapping cylinder structure 
$W\setminus Y\cong \partial W\times[0,1)$ we can further deform
the attaching map such that it maps to
$\partial W$.
Finally, using that $2i+1\le \dim(\partial W)=m-1$ we can  deform it to an embedding
(see \cite[Theorem II]{MR1503303}) into $\partial W$. We still denote this deformed attaching map by $\tilde \chi_\alpha$, and
we obtain a new deformed map $\chi_{\alpha}:=\tilde \chi_\alpha\times \{\alpha\}:S^{i}\to \partial W\times \{\alpha\}\subseteq \partial W\times \R$. 
Since the tangent bundle $T\partial W\to \partial W$ is trivial, we see that
the normal bundle $\nu_\alpha$ of the embedding $\tilde \chi_\alpha$ is stably trivial.
Since $\dim(\nu_\alpha)\ge i+1$ the bundle $\nu_\alpha$ is in fact trivial.

For each $\alpha\in J$ we now perform the procedure of attaching a handle to
$W\times \R$ described in  \cite[Sec .2]{MR0142128}. We can arrange the construction such that
for $\alpha\in J$ it takes place on $W\times (\alpha-1/4,\alpha+1/4)$.

The result of this construction is a manifold with boundary $(N,\partial N)$ of dimension
$m+1$ containing an $i+1$-dimensional $CW$-complex $Z$, and a map $N\to Z$  which represents 
$(N,\partial N)$  as a  tubular neighbourhood of $Z$, and  we have a homotopy equivalence
$Z\stackrel{\sim}{\to} \hat Z$. The construction depends on the choice of trivializations of the normal bundles $\nu_\alpha$. By Lemma \cite[Lemma 2.2]{MR0142128} we can choose the trivializations of the normal bundles $\nu_\alpha$ such that $TN\to N$ is again trivial.

After we finite iteration of this construction we obtain a manifold with boundary
$(\bar M,\partial \bar M)$ which is a tubular neighbourhood of a $CW$-complex $\tilde X$ which admits a homotopy equivalence $\tilde h:\tilde X\stackrel{\sim}{\to} X$.
The mapping cylinder structure on $\bar M$ gives rise to a projection $p:\bar M\to \tilde X$.
We now consider the smooth manifold $M:=\bar M\setminus \partial\bar M$.
The composition $\tilde h\circ p_{|M}:M\to X$
is a homotopy equivalence from a smooth manifold to $X$.
\hB

We can now construct the smooth $n$-connected approximation $j:\tilde EG\to EG$, where $n\ge 2$. 
We start with a countable $CW$-complex $BG$ of the homotopy type of the classifying space of $G$.
For example, we can take the standard simplicial model. It is countable since $G$ is countable.

We consider the $n$-skeleton $BG^n\subseteq BG$.
It is a finite-dimensional countable $CW$-complex. By Proposition \ref{nasowas} we can find a smooth manifold
$\tilde BG$ together with a homotopy equivalence 
$\hat j:\tilde BG\to BG^n$. Let $\bar j:\tilde BG\to BG$ denote the composition of $\hat j$ with the inclusion
$BG^n\hookrightarrow BG$. By construction $\bar j$ is $n$-connected.

Since $\bar j$ induces an isomorphism of fundamental groups it lifts to a $n$-connected map of universal coverings $j:\tilde EG\to EG$.

.

%
%
%

\subsection{The homological Chern character }

In this subsection we review the construction of the homological Chern
character given in \cite{MR1887884}, \cite{MR1914619}. 
Let $X$ be a proper $G$-CW-complex.
The main constituent of the Chern character 
is a homomorphism
$$\ch_H^X:H_{ev}(Z_G(H)\backslash X^H,\C)\otimes R_\C(H)\rightarrow K^G_0(X)$$
for any finite subgroup $H\subset G$.
\begin{eqnarray*}
H_{ev}(Z_G(H)\backslash X^H,\C)\otimes
R_\C(H)&\stackrel{(\pr_2)^{-1}_*\otimes \id}{\rightarrow}&
H_{ev}(EG\times_{Z_G(H)} X^H,\C)\otimes
R_\C(H)\\
&\stackrel{\ch^{-1}\otimes \id}{\rightarrow}&K_0(EG\times_{Z_G(H)}
X^H)_\C\otimes R_\C(H)\\
&\stackrel{\cong}{\rightarrow}&K_0^{Z_G(H)}(EG\times 
X^H)_\C\otimes K_0^H(*)_\C\\
&\stackrel{mult}{\rightarrow}&K_0^{Z_G(H)\times H}(EG\times 
X^H)_\C\\
&\stackrel{\Indu_{Z_G(H)\times H}^G}{\rightarrow}&K_0^G(\Indu^G_{Z_G(H)\times H}(EG\times 
X^H))_\C\\
&\stackrel{\Indu^G_{Z_G(H)\times H}(\pr_2)_*}{\rightarrow}&
K_0^G(\Indu^G_{Z_G(H)\times H}  X^H)_\C\\
&\stackrel{m_*}{\rightarrow}& K^G_0(X)
\end{eqnarray*}
Here $\ch$ is the homological Chern character, $\Indu^G_{Z_G(H)\times H}$ denotes the
induction functor, and $m:\Indu_{Z_G(H)\times H}^G  X^H=G\times_{Z_G(H)\times H}
X^H\rightarrow X$ is the $G$-map $(g,x)\mapsto gx$.

Let $C\subset G$ be a finite cyclic subgroup. Then we have a natural inclusion
$r^*:\C(\gen(C))\rightarrow R_\C(C)\cong \C C$ such that the image consists of
functions which vanish on $C\setminus \gen(C)$.
Note that $\C(\gen(C))$ and $H_*(Z_G(H)\backslash X^H,\C)$ are
left and right $W_G(C)$-modules in the natural way.
It follows from \cite{MR1914619}, Thm. 0.7,
that   
\begin{equation}\label{sur}\oplus_{(C)\in
C\cF Cyc(G)}\ch^X_C(1\otimes r^*):\bigoplus_{(C)\in
C\cF Cyc(G)}  H_{ev}(Z_G(C)\backslash
X^C,\C)\otimes_{\C W_G(C)} \C(\gen(C))\rightarrow K_0^G(X)_\C \end{equation}
is an isomorphism.

\section{Explicit decomposition of $K$-homology classes}

\subsection{An index formula}

Let $E$ be a $G$-equivariant vector bundle over $X$. If 
$A\in H_*(Z_G(H)\backslash
X^H,\C)\otimes  R_\C(H)$, then we can ask for a formula for
$\ind_\rho([E]\cap \ch^X_H(A))$ in terms of $\ch_X^H(b(\{E\}))$.
Let
$\epsilon:R(H)\rightarrow \Z$ be the homomorphism which takes the multiplicity
of the trivial representation. It extends to a group homomorphism
$\epsilon_\C:R_\C(H)\rightarrow \C$.  Using the ring structure of
$R_\C(H)$ and the pairing between homology and cohomology we obtain a
natural pairing 
\begin{eqnarray*}\langle.,.\rangle_\rho&:&\left(H_*(Z_G(H)\backslash
X^H,\C)\otimes  R_\C(H) \right)\otimes \left(H^*(Z_G(H)\backslash
X^H,\C)\otimes R_\C(H) \right) \\ &\rightarrow &
R_\C(H)\stackrel{\otimes [\rho_{|H}]}{\rightarrow}
R_\C(H)\stackrel{\epsilon_\C}{\rightarrow}  \C\ .\end{eqnarray*}

\begin{theorem}\label{inform}
$\ind_\rho([E]\cap \ch^X_H(A))=\langle\ch_X^H(b(\{E\})),A\rangle_\rho$
\end{theorem}
\proof
Let $M$ be a cocompact free even-dimensional  $Z_G(H)$-manifold equipped with
a invariant Riemannian metric and a Dirac operator $D$ associated to a
$Z_G(H)$-equivariant Dirac bundle $F\rightarrow M$. Furthermore, let
$f=(f_1,f_2):M\rightarrow EG\times X^H$ be a $Z_G(H)$-equivariant continuous
map. We form $[D]\in K_0^{Z_G(H)}(M)$ represented by the
Kasparov module $(L^2(M,F),\cF)$. Then $f_*[D]\in K_0^{Z_G(H)}(EG\times X^H)$.

Note that $K_0^{Z_G(H)}(EG\times X^H)_\C$ is spanned by elements arising in
this form. This can be seen as follows. First observe that every class in  
$K_0(EG\times_{Z_G(H)} X^H)$ can be represented in the form $\bar f_*[\bar D]$, where
 $\bar f:\bar N\to  EG\times_{Z_G(H)} X^H$ is a map from a closed $Spin^c$-manifold,
and $\bar D$ is the $Spin^c$-Dirac operator on $\bar N$. A proof of this result is given in 
\cite{math.KT/0701484}. We now consider
the pull-back
$$\xymatrix{N\ar[d]\ar[r]^f&EG\times X^H\ar[d]\\\bar N\ar[r]^{\bar f}&EG\times_{Z_G(H)} X^H}\ .$$
The manifold $N$ carries a $Z_G(H)$-invariant $Spin^c$-structure with associated Dirac operator $D$.
The class $f_*[D]$ corresponds to $\bar f_*[\bar D]$ under the isomorphism 
$K_0^{Z_G(H)}(EG\times X^H)\cong K_0^{Z_G(H)}(EG\times X^H)$.

Let $\phi\in \hat H$ be a finite-dimensional representation. It
gives rise to an element $[\phi]\in K_0^H(*)_\C$ under the natural
identification $R_\C(H)\cong K_0^H(*)_\C$. 
Let 
$T: K_0^{Z_G(H)}(EG\times X^H)\otimes K_0^H(*)_\C\rightarrow K_0^G(X)$ be the
composition $m_*\circ \Indu^G_{Z_G(H)\times H}(\pr_2)_*\circ 
\Indu^G_{Z_G(H)\times H}\circ mult$, which is part of the definition of
$\ch^X_H$. 

We first study $\ind_\rho([E]\cap T(f_*[D]\otimes[\phi]))$.
We have $mult\circ f_*([D]\otimes [\phi])=f_* \circ mult([D]\otimes [\phi])$,
and $mult([D]\otimes [\phi])\in K^{Z_G(H)\times H}_0(M)$ is represented by the
Kasparov module $(L^2(M,F)\otimes V_\phi,\cF\otimes \id)$.
Furthermore,
$\Indu_{Z_G(H)\times H}^G \circ f_* \circ mult([D]\otimes
[\phi])=\Indu_{Z_G(H)\times H}^G(f_*) \Indu_{Z_G(H)\times H}^G (
mult([D]\otimes [\phi]))$.
Explicitly, $\Indu_{Z_G(H)\times H}^G (
mult([D]\otimes [\phi]))$ is represented by a Kasparov module which is
constructed in the following way. Consider the exact sequence
$$0\rightarrow K\rightarrow Z_G(H)\times H\rightarrow G\ ,$$
where $K=Z_G(H)\cap H= Z_H(H)$.  We identify $K\backslash Z_G(H)\times H$ with the subgroup
$Z_G(H)H\subseteq G$.

Note that we consider $M$ as a $Z_G(H)\times H$-manifold via the action of the first factor.
The $Z_G(H)H$-manifold $\hat M:=K\backslash M$ carries an induced equivariant 
Dirac bundle $\hat F$. We further consider the flat $Z_G(H)H$-equivariant 
bundle $\hat V_\phi:= V_\phi\times_K M$ over $\hat M$. The twisted bundle $\hat
F\otimes \hat V_\phi $ is a $Z_G(H)H$-equivariant Dirac bundle.
We consider the cocompact proper $G$-manifold $\tilde M :=G\times_{Z_G(H)H}
\hat M$. The $Z_G(H)H$-equivariant Dirac bundle $\hat F\otimes \hat V_\phi$
induces a $G$-equivariant Dirac bundle $\tilde F_\phi\rightarrow \tilde M$  in
a natural way with associated operator $\tilde D_\phi$. Then
$\Indu_{Z_G(H)\times H}^G ( mult([D]\otimes [\phi]))$ is represented by
$[\tilde D_\phi]$. The map $\Indu_{Z_G(H)\times H}^G(f_*)$ is induced by the
$G$-map $\tilde f :\tilde M\rightarrow G\times_{Z_G(H)H} (K\backslash
EG\times  X^H)$ given by $\tilde f([g,Km]):=[g , (Kf_1(m),f_2(m))]$. It is now
clear that $T(f_*[D]\otimes[\phi])$ is represented by $h_*[\tilde D_\phi]$,
where $h:\tilde M\rightarrow X$ is given by $h([g,Km])=g f_2(m)$.

It follows from the associativity of the Kasparov product that 
$$\ind_\rho([E]\cap T(f_*[D]\otimes[\phi]))=\ind_\rho([E] \cap
h_*[\tilde D_\phi])=\ind_\rho([h^* E]\cap [\tilde D_\phi])\ .$$
By Theorem \ref{orbi} and Proposition \ref{cca}
we obtain 
$\ind_\rho([h^* E]\cap [\tilde D_\phi])=\ind(\bar D_{\phi,h^*E,\rho})$,
where $\tilde D_{\phi,h^*E}$ is the $G$-invariant  Dirac operator associated
to $\tilde F\otimes h^*E$, and $\bar D_{\phi,h^*E,\rho}$ is the  operator on
the orbifold $\bar M:=G\backslash \tilde M$ induced by $\tilde D_{\phi,h^*E}$
and the twist $\rho$.  Restriction from $\tilde M$ to the submanifold
$\{1\}\times \hat M$   provides an isomorphism
$$\left(C^\infty(\tilde M,\tilde
F_\phi\otimes h^*E)\otimes V_\rho\right)^G\cong \left(C^\infty(\hat M,\hat F
\otimes \hat V_\phi\otimes \bar f_2^*E_{|X^H})\otimes
V_\rho\right)^{Z_G(H)\times H}\ ,$$ where $\bar f_2:\hat M\rightarrow X^H$ is
induced by $f_2$. 
Since the action of $H$ on the latter
spaces is implemented by the action on the fibres of $V_\phi\otimes
\bar f_2^*E_{|X^H}\otimes V_\rho$ we further obtain $$\left(C^\infty(\tilde
M,\tilde F_\phi\otimes h^*E)\otimes V_\rho\right)^G= C^\infty(\hat M, \hat F
\otimes (V_\phi\otimes \bar f_2^*E_{|X^H}\otimes V_\rho)^H)^{K\backslash
Z_G(H)}\ .$$  
In the present situation we
have $\bar M=Z_G(H)\backslash M=(K\backslash Z_G(H))\backslash \hat M$, i.e.
the orbifold is smooth, and it carries the Dirac bundle $\bar F$ with
associated Dirac operator $\bar D$. We define the
$ (K\backslash Z_G(H))$-equivariant bundle $E_{\phi\otimes\rho}:=(V_\phi\otimes
E_{|X^H}\otimes V_\rho)^H$ over $X^H$. Furthermore, we consider the quotient 
$\overline{\bar f_2^* E_{\phi\otimes\rho}}:=(K\backslash Z_G(H))\backslash
\bar f_2^* E_{\phi\otimes\rho}$ over $\bar M$. The identifications above show
that $\ind_\rho(\bar D_{\phi,h^*E})=\ind(\bar D_{\overline{\bar f_2^* E_{\phi\otimes\rho}}})$,
i.e. it is the index of a twisted Dirac operator. Writing the index of the
twisted Dirac operator in terms of Chern characters we obtain
$$\ind_\rho([E]\cap T(f_*[D]\otimes[\phi]))=\langle \ch(\overline{ \bar f_2^*
E_{\phi\otimes\rho}}),\ch([\bar D])\rangle \ .$$  Note that  
$\overline{\bar f_2^* E_{\phi\otimes\rho}}=\bar f^* \overline{\pr_2^*
E_{\phi\otimes\rho}} $, where $\pr_2:EG\times X^H\rightarrow X^H$, $\bar
f:\bar M\rightarrow  EG\times_{Z_G(H)} X^H$ is induced by $f$, and 
$  \overline{\pr_2^* E_{\phi\otimes\rho}}:=Z_G(H)\backslash \pr_2^*
E_{\phi\otimes\rho}$. We conclude
that $$\langle \ch(\overline{ \bar f_2^*
E_{\phi\otimes\rho}}),\ch([\bar D])\rangle
=\langle \ch(\overline{ \pr_2^* E_{\phi\otimes\rho}}),
 \ch(\bar f_*[\bar D])\rangle\ .$$
The right-hand side can now be written as
$$\left\langle \epsilon_\C\left(\ch^H_X(b([E]))\otimes[\phi]\otimes\rho\right),
 (\pr_2)_*\ch(\bar f_*[\bar D])\right\rangle=\langle \ch^H_X(E),(\pr_2)_* \ch(\bar f_*[\bar
D])\otimes [\phi])\rangle_\rho \ .$$
Note that $\ch_H^X((\pr_2)_*\ch(\bar f_*[\bar D])\otimes
[\phi])=T(f_*[D]\otimes [\phi])$.
Therefore we have shown
$$\ind_\rho([E]\cap  \ch_H^X((\pr_2)_*\ch(\bar f_*[\bar D])\otimes
[\phi])))=  \langle \ch^H_X(b([E])),\ch(\bar f_*[\bar D])\otimes
[\phi] \rangle_\rho \ .$$
Since the classes  $(\pr_2)_*\ch(\bar f_*[\bar D])\otimes [\phi]$ for varying data
$M,F,f,\phi$ span $H_{ev}(Z_G(H)\backslash X^H,\C)\otimes R_\C(H)$
the theorem follows.
\hB

\subsection{Decomposition}

\begin{lem}\label{non}
Let $X$ be a finite proper $G$-CW-complex.
If $x\in K^G_0(X)_\C$ and
$\ind([E]\cap x)=0$ for all $G$-equivariant complex vector bundles $E$ on $X$,
then $x=0$.
\end{lem}
\proof
Because of the isomorphism (\ref{sur}) it suffices to show that if $A\in
H_{ev}(Z_G(C)\backslash X^C,\C)\otimes_{\C W_G(C)} \C(\gen(C))$ and
$\ind([E]\cap \ch_C^X(A))=0$ for all $E$, then $A=0$.
By Theorem \ref{inform} we have
$\ind([E]\cap \ch_C^X(A))=\langle \ch^C_H(b(\{E\})),A\rangle$.
Using the surjectivity of $b$ and of the isomorphism (\ref{sur1}), and the fact
that the pairing $$\langle.,.\rangle :\left(H^{ev}(Z_G(C)\backslash
X^C,\C)\otimes \C(\gen(C))\right)^{W_G(C)}\otimes\left(H_{ev}(Z_G(C)\backslash
X^C,\C)\otimes_{\C W_G(C)} \C(\gen(C))\right)$$
is nondegenerate we see that 
$\langle \ch^C_H(b(\{E\})),A\rangle=0$  for all $E$ indeed implies $A=0$. \hB

Let now $M$ be an even-dimensional proper cocompact $G$-manifold equipped
with a $G$-invariant Riemannian metric $g^M$ and a $G$-equivariant Dirac
bundle $F$ with associated Dirac operator $D$. Let $[D]_\C\in K_0^G(M)_\C$ be
the equivariant $K$-homology class of $D$.

The $G$-space $M$ has the $G$-homotopy type of a
finite proper $G$-CW-complex. In particular, we have
the isomorphism (\ref{sur}) 
$$\oplus_{(C)\in
C\cF Cyc(G)}\ch^M_C(1\otimes r^*):\bigoplus_{(C)\in
C\cF Cyc(G)}  H_{ev}(Z_G(C)\backslash
M^C,\C)\otimes_{\C W_G(C)} \C(\gen(C))\rightarrow K_0^G(M)_\C \ .$$
Therefore, there exist uniquely determined  classes $[D](C)\in
H_{ev}(Z_G(C)\backslash M^C,\C)\otimes_{\C W_G(C)} \C(\gen(C))$ such that
$$\sum_{{(C)\in
C\cF Cyc(G)}} \ch^M_C(1\otimes r^*)([D](C))=[D]_\C\ .$$
\begin{theorem}\label{hhjj}
We have the equality 
$$[D](C)= [\hat U]\ ,$$ 
where $[\hat U]$ is given by $ \gen(C)\ni g
\rightarrow [\hat U(g)]\in H_{ev}(Z_G(C)\backslash M^C,\C)$,
and $\hat U(g)$ was defined in (\ref{hatudef}). 
\end{theorem}
\proof
Let $E$ be any $G$-equivariant complex vector bundle over $M$.
Then we have
$$\ind([E]\cap [D]_\C)=\sum_{{(C)\in
C\cF Cyc(G)}} \langle \ch_M^C(b(\{E\})),[D](C)\rangle\ .$$
Using the definition of $\epsilon_\C$ and Proposition \ref{chern}
we can write out the summands of right-hand side as follows
$$\langle \ch_M^C(b(\{E\})),[D](C)\rangle = \frac{1}{|C|}\sum_{g\in \gen(C)}
\langle [\ch(g,E)],[D](C)(g)\rangle \ .$$
On the other hand the index formula Theorem \ref{oopp}
gives
$$\ind([E]\cap [D]_\C)=\sum_{{(C)\in
C\cF Cyc(G)}}\frac{1}{|C|} \sum_{g\in \gen(C)}
 \langle [\ch(g,E)],[\hat U](g) \rangle\ .$$
Varying $E$ and using Lemma \ref{non} we conclude
$[D](C)=[\hat U]$.
\hB


\begin{thebibliography}{GHT00}

\bibitem[BGV92]{MR1215720}
Nicole Berline, Ezra Getzler, and Mich{\`e}le Vergne.
\newblock {\em Heat kernels and {D}irac operators}, volume 298 of {\em
  Grundlehren der Mathematischen Wissenschaften [Fundamental Principles of
  Mathematical Sciences]}.
\newblock Springer-Verlag, Berlin, 1992.

\bibitem[BHS]{math.KT/0701484}
Paul Baum, Nigel Higson, and Thomas Schick.
\newblock {On the Equivalence of Geometric and Analytic K-Homology},
  arXiv:math.KT/0701484.


\bibitem[Bla98]{MR1656031}
Bruce Blackadar.
\newblock {\em {$K$}-theory for operator algebras}, volume~5 of {\em
  Mathematical Sciences Research Institute Publications}.
\newblock Cambridge University Press, Cambridge, second edition, 1998.

\bibitem[BM04]{MR2101228}
Paul Balmer and Michel Matthey.
\newblock Model theoretic reformulation of the {B}aum-{C}onnes and
  {F}arrell-{J}ones conjectures.
\newblock {\em Adv. Math.}, 189(2):495--500, 2004.

\bibitem[Bre93]{MR1224675}
Glen~E. Bredon.
\newblock {\em Topology and geometry}, volume 139 of {\em Graduate Texts in
  Mathematics}.
\newblock Springer-Verlag, New York, 1993.

\bibitem[Bro62]{MR0142128}
Edgar~H. Brown, Jr.
\newblock Nonexistence of low dimension relations between {S}tiefel-{W}hitney
  classes.
\newblock {\em Trans. Amer. Math. Soc.}, 104:374--382, 1962.



\bibitem[Bun95]{MR1348799}
Ulrich Bunke.
\newblock A {$K$}-theoretic relative index theorem and {C}allias-type {D}irac
  operators.
\newblock {\em Math. Ann.}, 303(2):241--279, 1995.

\bibitem[Che73]{MR0369890}
Paul~R. Chernoff.
\newblock Essential self-adjointness of powers of generators of hyperbolic
  equations.
\newblock {\em J. Functional Analysis}, 12:401--414, 1973.

\bibitem[DL98]{MR1659969}
James~F. Davis and Wolfgang L{\"u}ck.
\newblock Spaces over a category and assembly maps in isomorphism conjectures
  in {$K$}- and {$L$}-theory.
\newblock {\em $K$-Theory}, 15(3):201--252, 1998.

\bibitem[Far92a]{MR1127139}
Carla Farsi.
\newblock {$K$}-theoretical index theorems for good orbifolds.
\newblock {\em Proc. Amer. Math. Soc.}, 115(3):769--773, 1992.

\bibitem[Far92b]{MR1164622}
Carla Farsi.
\newblock {$K$}-theoretical index theorems for orbifolds.
\newblock {\em Quart. J. Math. Oxford Ser. (2)}, 43(170):183--200, 1992.

\bibitem[Far92c]{MR1163555}
Carla Farsi.
\newblock A note on {$K$}-theoretical index theorems for orbifolds.
\newblock {\em Proc. Roy. Soc. London Ser. A}, 437(1900):429--431, 1992.

\bibitem[Fed88]{0780.54034}
V.V. Fedorchuk.
\newblock {The fundamentals of dimension theory.}
\newblock {General topology. I. Basic concepts and constructions. Dimension
  theory. Encycl. Math. Sci. 17, 91-192 (1990); translation from Itogi Nauki
  Tekh., Ser. Sovrem. Probl. Mat., Fundam. Napravleniya 17, 111-224 (1988).},
  1988.


\bibitem[GHT00]{MR1711324}
Erik Guentner, Nigel Higson, and Jody Trout.
\newblock Equivariant {$E$}-theory for {$C\sp *$}-algebras.
\newblock {\em Mem. Amer. Math. Soc.}, 148(703):viii+86, 2000.

\bibitem[Kas88]{MR918241}
G.~G. Kasparov.
\newblock Equivariant {$KK$}-theory and the {N}ovikov conjecture.
\newblock {\em Invent. Math.}, 91(1):147--201, 1988.

\bibitem[Kaw78]{MR0474432}
Tetsuro Kawasaki.
\newblock The signature theorem for {$V$}-manifolds.
\newblock {\em Topology}, 17(1):75--83, 1978.

\bibitem[Kaw79]{MR527023}
Tetsuro Kawasaki.
\newblock The {R}iemann-{R}och theorem for complex {$V$}-manifolds.
\newblock {\em Osaka J. Math.}, 16(1):151--159, 1979.

\bibitem[Kaw81]{MR641150}
Tetsuro Kawasaki.
\newblock The index of elliptic operators over {$V$}-manifolds.
\newblock {\em Nagoya Math. J.}, 84:135--157, 1981.

\bibitem[LO01a]{MR1851256}
Wolfgang L{\"u}ck and Bob Oliver.
\newblock Chern characters for the equivariant {$K$}-theory of proper
  {$G$}-{CW}-complexes.
\newblock In {\em Cohomological methods in homotopy theory (Bellaterra, 1998)},
  volume 196 of {\em Progr. Math.}, pages 217--247. Birkh\"auser, Basel, 2001.

\bibitem[LO01b]{MR1838997}
Wolfgang L{\"u}ck and Bob Oliver.
\newblock The completion theorem in {$K$}-theory for proper actions of a
  discrete group.
\newblock {\em Topology}, 40(3):585--616, 2001.

\bibitem[L{\"u}c02a]{MR1887884}
Wolfgang L{\"u}ck.
\newblock Chern characters for proper equivariant homology theories and
  applications to {$K$}- and {$L$}-theory.
\newblock {\em J. Reine Angew. Math.}, 543:193--234, 2002.

\bibitem[L{\"u}c02b]{MR1914619}
Wolfgang L{\"u}ck.
\newblock The relation between the {B}aum-{C}onnes conjecture and the trace
  conjecture.
\newblock {\em Invent. Math.}, 149(1):123--152, 2002.

\bibitem[MN06]{MR2193334}
Ralf Meyer and Ryszard Nest.
\newblock The {B}aum-{C}onnes conjecture via localisation of categories.
\newblock {\em Topology}, 45(2):209--259, 2006.


\bibitem[Whi36]{MR1503303}
Hassler Whitney.
\newblock Differentiable manifolds.
\newblock {\em Ann. of Math. (2)}, 37(3):645--680, 1936.

\bibitem[Whi55]{MR0073980}
Hassler Whitney.
\newblock On singularities of mappings of euclidean spaces. {I}. {M}appings of
  the plane into the plane.
\newblock {\em Ann. of Math. (2)}, 62:374--410, 1955.

\end{thebibliography}
\end{document}